\documentclass[11
pt]{article}
\usepackage[T1]{fontenc}
\usepackage[english]{babel}

\usepackage{latexsym}
\usepackage{amssymb}
\usepackage{stmaryrd}
\usepackage{amsmath}
\usepackage{paralist} %liste numerotee
\usepackage[dvips]{graphicx,epsfig}
\usepackage{color}

\newtheorem{defi}{Definition}

\newtheorem{thm}{Theorem}

\newtheorem{cor}{Corollary}
                      
\newtheorem{rem}{Remark}
\newcommand{\conv}[1]{\textrm{conv}}
\newcommand{\var}[1]{\textrm{Var}}
\newcommand{\CQFD}{\nolinebreak \hfill\rule{2mm}{2mm}\medbreak\par}                                                                               
\newcommand{\vect}[1]{\textrm{vect}}
\newcommand{\supp}[1]{\textrm{supp}}
\newcommand{\sth}[1]{\textrm{such that}}

\begin{document}
\begin{center}
{\Large \textbf{On convex hull and winding number of 
self similar processes}}\\
\vspace{0.5cm}

{\large Yu. Davydov\footnote{University  Lille 1, Laboratory P. Painlev\'e, France}
}
\end{center}

\vspace{1.5cm}

\noindent\textbf{Abstract:} {\small
\vspace{5pt}

It is well known that  for a standard Brownian motion (BM) $ \{B(t), \;t \geq 0\}$ with values in $\mathbb{R}^d$,
its convex hull $ V(t)=\conv \{\{\,B(s),\;s \leq t \}$ with probability $1$ for each $t > 0$ contains $0$ as  an interior point  (see Evans (1985)).
We also know that the winding number of a typical path of a $2$-dimensional BM is equal to $+\infty.$

 The aim of this article is to show that these properties aren't specifically "Brownian", but hold
 %state the analogical properties 
 for a much larger class of $d$-dimensional self similar processes. This class contains in particular
$d$-dimensional fractional Brownian motions and (concerning convex hulls) strictly stable Levy processes.    
}

\vspace{10pt}

\vspace{5pt}

%------------------------------------------------------------------------------
\section{Introduction}
%------------------------------------------------------------------------------

Let $(\Omega, \cal{F}, \mathbb{P})$ be a basic probability space. Consider a $d$-dimensional  process 
$X = \{X(t),\; t\geq 0\}$ defined on $\Omega$ which is self-similar of index $H>0.$\\
It means that for each constant $c>0$ the process 
$\{X(ct),\; t\geq 0\}$ has the same distribution as 
$\{c^HX(t),\; t\geq 0\}.$ 
\vspace{5pt}

Let $L = \{L(u),\; u\in \mathbb{R}^1\}$ be the strictly stationary process obtained from $X$ by Lamperti transformation:

\begin{equation}\label{L}
\hspace{30pt} L(u) = e^{-Hu} X(e^u),\;\;\; u \in \mathbb{R}^1.
\end{equation}

Equivalently,

$$
\hspace{30pt} X(t)=t^HL(\log {t}),\;\;\; t \in \mathbb{R}_*^+.
$$

Let $\Theta = \{0,1\}^d$ be the set of all dyadic sequences of length $d.$
Denote by $D_\theta,$\\$ \;\theta \in \Theta,$ the quadrant 
$$
D_\theta = \prod_{i=1}^{d}\mathbb{R}_{\theta_i},
$$
where $\mathbb{R}_{\theta_i} = [0,\infty)$ if  $\theta_i = 1,$ and $\mathbb{R}_{\theta_i} = (-\infty, 0]$ if
$\theta_i = 0.$
\vspace{5pt}

The positive quadrant $D_{(1,1,\ldots,1)}$ for simplicity is denoted by $D.$	
\vspace{5pt}

We say that the process $X$ is {\it non degenerate} if
for all $\theta \in \Theta$	
$$
\mathbb{P}\{X(1)\in D_\theta \}>0.
$$  

Two important examples of self similar processes are \textbf{fractional Brownian motion}
and \textbf{stable Levy process}.
%\vspace{5pt}
	
\begin{defi}
We call a self-similar (of index $H>0$) process $B^H$
{\it fractional Brownian motion} (FBM) if for each $e \in \mathbb{R}^d$ the scalar process 
$t \rightarrow \langle B^H(t),e\rangle$ is a standard one-dimensional FBM of index $H$ up to a constant $c(e).$
\end{defi}
\vspace{5pt}

It is easy to see that in this case $c^2(e) = \langle Qe,e\rangle,$ where $Q$ is the covariance matrix of $B^H(1),$ and hence
%\begin{block}{}
$$
\mathbb{E}\langle B^H(t),e\rangle\langle B^H(s),e\rangle = \langle Qe,e\rangle\frac{1}{2}(t^{2H} +s^{2H}-|t-s|^{2H}),\;\;\; t,s \geq 0;\,\, e\in \mathbb{R}^d.
$$ 
The process $B^H$ is non degenerate iff the rank of the matrix $Q$ is equal to $d$.
If $H=\frac{1}{2}, Q = I_d,$ then $B^H$ is a standard Brownian motion.
\vspace{5pt}
	
(See A. Xiao (2013), Ra$\breve c$kauskas and Ch. Suquet (2011), F. Lavancier et al. (2009) and references therein for more general definitions of operator self-similar FBM).
%\vspace{5pt}
	
\begin{defi}
We call $S= \{S(t), t\in \mathbb{R}_+\}$ {\it $\alpha$-strictly stable Levy process} (StS) if
\begin{itemize}
\item[1)] $S(1)$ has a $\alpha$-strictly stable distribution in $\mathbb{R}^d$;
\item[2)] it has independent and stationary increments;
\item[3)] it is continuous in probability.
\end{itemize}
\end{defi}
Then for each $t\in \mathbb{R}_+$ the random variable $S(t)$ has the same distribution as 
$t^{\frac{1}{\alpha}}S(1).$ 

The cadlag version of $S$ on $[0,1]$ can be obtained with the help of LePage series representation (see \cite{ST} for more details).
If $\alpha \in (0, 1)$ or if $\alpha \in (1, 2)$ and $EX(1) = 0,$
then we have:
\begin{equation}\label{LePage}
\{S(t), t\in [0,1]\} \stackrel{\cal L}{=} 
\{c\sum_1^\infty \Gamma_k^{-1/\alpha}\varepsilon_k{\mathbf 1}_{[0,t]}
(\eta_k),\;\;t\in [0,1]\},
\end{equation}
where $c$ is a constant, $\Gamma_k = \sum_1^k\gamma_j,$ 
$\{\gamma_j\}$ is a sequence of i.i.d. random variables with common standard exponential distribution, $\{\varepsilon_k\}$ is 
a sequence of i.i.d. random variables with common distribution
$\sigma$ concentrated on unit sphere $\mathbb{S}^{d-1}$,
$\{\eta_k\}$ is a sequence of $[0,1]$-uniformly distributed i.i.d. random variables, and the three sequences $\{\gamma_j\},\;\{\varepsilon_k\},\;\{\eta_k\}$ are supposed to be independent.

The measure $\sigma$ is called spectral measure of $S$.
It is easy to see that if (\ref{LePage}) takes place, the process 
$X$ is non degenerate iff 
$\;\;\textrm{vect} \{ \textrm{supp}\; {\sigma}\}= \mathbb{R}^1.$
\vspace{7pt}

In Section 2 the object of our interest is the convex hull process
$V = \{V(t)\}$ associated with $X.$ We show that under very sharp
conditions with probability 1 for all $t>0$ the convex set $V(t)$ contains $0$ as its interior point. From this result some interesting corollaries are deduced.

Section 3 is devoted to studying the winding numbers of two-dimensional self similar processes. As a corollary of our main result we show that for the typical path of a standard two-dimensional FBM the number of its clockwise and anti-clockwise winds around $0$ in the neighborhood of zero or at infinity is equal to
$+\infty.$
\section{Convex hulls}

For a 
Borel set $A \subset \mathbb{R}^d$ we denote by $\mathrm{conv}{(A)} $ the closed convex hull of $A$ and define
%In this section the object of our interest is
 the convex hull process related to $X$:
%\begin{equation}
$$
V(t) = \mathrm{conv} \{\,X(s),\;\;s\leq t\}.
$$
%\end{equation}   

\begin{thm}\label{th1} 
 \it{Let $X$ be a non degenerate self similar process such that 
  the strictly stationary process $L$ generating $X$ is ergodic. 
 Then with probability $1$ for all $t>0$ the point $0$ is an interior point of $V(t).$}
\end{thm}

{\bf Application to FBM.}
Let $B^H$ be a FBM with index $H.$
The next  properties follow from the definition without difficulties.

\begin{enumerate}
\item[1)] {\bf Continuity.} The process $X$ has a continuous version.

\noindent
Below we always suppose $B^H$ to be continuous.

\item[2)] {\bf Reversibility.} If the process $Y$ is defined by
$$
Y(t) = B^H(1) -B^H(1-t),\;\;\; t\in [0,1],
$$
then $\{Y(t),\; t\in [0,1]\} \stackrel{\cal L}{=} \{B^H(t),\; t\in [0,1]\},$ 
where $\stackrel{\cal L}{=}$ means equality in law.

\item[3)] {\bf Ergodicity.} Let $L = \{L(u),\; u\in \mathbb{R}^1\}$ be the strictly stationary Gaussian process obtained from $B^H$ by Lamperti transformation (\ref{L}).

Then $L$ is ergodic (see Cornfeld et al. (1982), Ch. 14, \S 2,  Th.1, Th.2).
\end{enumerate}		

It is supposed below  that the process  $B^H$ is non degenerate.

\begin{cor} \label{cor1} 
 \it{Let $V$ be the convex hull process related to $B^H.$ Then with probability $1$ for all $t>0$ the point $0$ is an interior point of $V(t).$}
\end{cor}

This follows immediately from Th.1.

\begin{cor} \label{cor2} 
 \it{Let $V$ be the convex hull process related to $B^H.$ Then for each $t>0$ with probability $1$  the point $B^H(t)$ is an interior point of $V(t).$}
\end{cor}

{\bf Proof of Corollary 2.} Denote by $A^\circ $ the interior of $A.$ By self-similarity of the process $B^H$ it is sufficient to state this property for $t=1.$ Then, due to the reversibility of $B^H$ by Th. 1., a.s. 
\begin{equation}
\label{revers}
0 \in [\mathrm{conv}\{\,B^H(1) -B^H(1-t),\;\;\; t\in [0,1]\}]^\circ.
\end{equation}

As
$$
\mathrm{conv}\{\,B^H(1) -B^H(1-t),\;\; t\in [0,1]\} = B^H(1) -\conv\{\{\,B^H(1-s),\;\;\; s\in [0,1]\},
$$
the relation  (\ref{revers}) is equivalent to
$$
B^H(1) \in [\mathrm{conv}\{\,B^H(s),\;\;s\in [0,1]\}]^\circ,
$$
 which concludes the proof.\CQFD

%%%%%%%%%%%%%%%%%%%%%%%%%%%%%%%%%%%%%%%%%%%%%%%%%%%%%%%%

Let ${\cal K}_d$ be the family of all compact convex subsets of $\mathbb{R}^d.$ It is well known that
${\cal K}_d$ equipped with Hausdorff metric is a Polish space.
\vspace{5pt}

We say that a function $f : [0,1] \rightarrow {\cal K}_d$ is {\it increasing}, if $f(t) \subset f(s)$ for
$0\leq t<s\leq 1.$
\vspace{5pt}

We say that a function $f : [0,1] \rightarrow {\cal K}_d$ is {\it almost everywhere constant},  
if $f$ is  such that for almost every $t\in [0,1]$ there exists an interval $(t-\varepsilon, t+\varepsilon)$ where $f$ is constant.
\vspace{5pt}

We say that a function $f : [0,1] \rightarrow {\cal K}_d$ is a {\it Cantor - staircase} (C-S), 
%(or Devil's staircase), 
if $f$ is continuous, increasing and almost everywhere constant.
\vspace{5pt}

The next statement  follows easily  from Corollary 2.

\begin{cor}  \label{cor3} 
 \it{Let $V$ be the convex hull process related to $B^H.$ Then with probability $1$ the paths of the process $t\rightarrow V(t)$ are C-S functions.}
 \end{cor}
 
 \begin{rem} \label{th1-2}
Let $h :{\cal K} \rightarrow \mathbb{R}^1 $ be an increasing continuous function. Then  almost all paths of the process $t \rightarrow h(V(t))$ are C-S real  functions. This obvious fact may be applied to all reasonable geometrical characteristics of $V(t),$ such as volume, surface area, diameter,...
\end{rem}
\vspace{5pt}

{\bf Application to StS.} 
Let now $S$ be a StS process with exponent $\alpha <2.$
The following  properties are more or less known.

\begin{enumerate}
\item[1)] {\bf Right continuity.} The process $S$ has a {\it cadlag} version (see remark above just after the definition).

\vspace{7pt}

\item[2)] {\bf Reversibility.} Let
$$
Y(t) = S(1) -S(1-t),\;\;\; t\in [0,1].
$$
Then $\{Y(t),\; t\in [0,1]\} \stackrel{\cal L}{=} \{S(t),\; t\in [0,1]\}.$ 

\vspace{7pt}

\item[3)] {\bf Self-similarity.} The process $S$ is self-similar of index $H = \frac{1}{\alpha}.$
\vspace{7pt}

\item[4)] {\bf Ergodicity.} Let $L = \{L(u),\; u\in \mathbb{R}^1\}$ be the strictly stationary  process obtained from $S$ by Lamperti transformation (\ref{L}).
Then $L$ is ergodic.
\end{enumerate}

 We  suppose that the law of $S(1)$  is non degenerate.
\begin{cor}  \label{cor4} 
 \it{Let $V$ be the convex hull process related to $S.$  Then with probability $1$ for all $t>0$ the point $0$ is an interior point of $V(t).$}
\end{cor} 

\begin{cor} \label{cor5} 
 \it{Let $V$ be the convex hull process related to $S.$ Then for each $t>0$ with probability $1$  the point $X(t)$ is an interior point of $V(t).$}
\end{cor}

\begin{cor}  \label{cor6} 
 \it{Let $V$ be the convex hull process related to $S.$ Then with probability $1$ the paths of the process $t\rightarrow V(t)$ are right continuous almost everywhere constant functions.}
\end{cor}
We  omit proofs of these statements as they are similar to proofs of Corollaries 1 - 3.
\vspace{7pt}

{\bf Proof of Theorem 1.}
We  first show that
\begin{equation}
p\;\; \stackrel{def}{=}\;\; \mathbb{P}\{\,\exists\, t\in (0,1]\;\vert\; X(t) \in D^\circ\} = 1.
\end{equation}
Remark that $p$ is strictly positive:
\begin{equation}
\label{p}
p\geq \mathbb{P}\{X(1) \in D^\circ\} > 0
\end{equation}
due to the hypothesis that the law of $X(1)$ is non degenerate.

By self similarity
$$
\mathbb{P}\{D^\circ\ \cap \{X(t), t \in [0,T]\}= \emptyset\} = 1-p
$$
for every $T>0.$
The sequence of events $(A_n)_{n\in\mathbb{N}}\,,$
$$
A_n = \{D^\circ\ \cap \{X(t), t \in [0,n]\}= \emptyset\},
$$
being decreasing, it follows that
$$
1-p = \lim \mathbb{P}(A_n) = \mathbb{P}(\cap_nA_n) = \mathbb{P}\{X(t) \notin D^\circ,\;\; \forall t\geq0\}.
$$
In terms of the stationary process $L$ from Lamperti representation %(\ref{lamperti}) 
it means that
$$
\mathbb{P}\{L(s) \notin D^\circ,\;\; \forall s \in \mathbb{R}^1 \}= 1-p.
$$
As this event is invariant, by ergodicity of $L$ and due to (\ref{p}) we see that the value $p=1$ is the  only one possible.

Applying the similar arguments to another  quadrants 
$D_\theta, \theta \in \Theta$, we get that with probability 1 there exists points  $t_\theta \in (0,1],$ such that\\ $X(t_\theta) \in D_\theta^\circ, \;\;\; \theta \in \Theta.$
Now, to end the proof it is sufficient to remark that
$$
V(1)^\circ = \mathrm{conv}\{X(t), t \in [0,1]\}^\circ \supset \mathrm{conv} \{X(t_\theta), \theta \in \Theta\}^\circ
$$
and that the last set evidently contains $0.$ \CQFD

\section{Winding numbers}
Nown we consider a  2-dimensional self similar process
$X = \{X(t),\; t\geq 0\}.$  It is supposed that the following properties are fulfilled:
\begin{enumerate}
\item[1)] Process $X$ is continuous.
\item[2)] Process $X$ is non degenerate.
\item[3)] Process $X$ is symmetric: $X$ and $-X$ have the same law.
\item[4)] The stationary process $L$ associated with $X$ is ergodic.
\item[5)] Starting from $0$ the process $X$  with probability 1 never come back:
\begin{equation}\label{nonzero}
\mathbb{P}\{X(t) \neq 0,\;\; \forall\; t>0 \}= 1.
\end{equation}
\end{enumerate}
Due to the last hypothesis, considering  $\mathbb{R}^2$
as complex plane, we can define the winding numbers (around $0$)
$\;\;\nu[s,t], \,0<s<t,$ by the usual way (see \cite{Y}, Ch.5): 
%as the increments of the $\arg{(B^H(t)}:$
$$
\nu[s,t] = \arg{(X(t))} - \arg{(X(s))}.
$$
We set
$$
\nu_+(0,t] = \limsup_{s\downarrow 0}\nu[s,t],\;\;\;
\nu_-(0,t] = \liminf_{s\downarrow 0}\nu[s,t],
$$
$$
\nu_+[s,\infty) = \limsup_{t\rightarrow \infty}\nu[s,t],\;\;\;
\nu_-[s,\infty) = \liminf_{t\rightarrow \infty}\nu[s,t].
$$
The values $\nu_+(0,t],\;\nu_-(0,t]$ represent respectively the number of clockwise and anti-clockwise winds around $0$ in the neighborhood of the starting point, while $\nu_+[s,\infty),\;
\nu_-[s,\infty)$ are the similar winding numbers at infinity.

\begin{thm}\label{theorem2}
Let $X$ be a  2-dimensional self similar process with the properties 1) -- 5) mentioned above. 
Then with probability one for all $t>0$
\begin{equation}\label{winding}
\nu_+(0,t] = \nu_+[t,\infty)= -\nu_-(0,t] = - \nu_-[t,\infty)  = + \infty.
\end{equation}
\end{thm} 
\begin{cor}\label{fbm-wind}
Let $B^H$ be a  2-dimensional standard FBM and assume that $H\in [1/2,\,1).$ Then with probability one for all $t>0$ the equalities (\ref{winding}) take place.
\end{cor}
{\bf Proof.} Case $H=1/2$ is well known, see \cite{Y}, Ch. 5, which contains exhaustive information on Brownian winding numbers.

If $H\in (1/2,\,1),$ we apply  Theorem 2 as all
hypothesis 1)--5) are fulfilled: indeed, the properties 1)--3) are obvious; the ergodicity of $L,\; L(t) = (L_1(t), L_2(t)),$ follows from the fact that $EL_1(t)L_1(0)\rightarrow 0$ when $t \rightarrow \infty$ (see, \cite{CFS}, Ch. 14, Sec. 2, Th.2); The property 5) can be deduced from Th. 11 of \cite{X} (see also Th. 4.2 of \cite{X1} and
Th. 2.6 of \cite{X2}).\CQFD

\begin{rem}
If $H \in (0, \frac{1}{2}),$ the process $t\rightarrow 
\arg B^H(t)-\arg B^H(0)$ is not continuous with positive probability
as the set $\{t\in (0,1]\;|\; B^H(t) = 0\}$ is not empty 
(see \cite{X},  Th. 11). It means that in this case the winding numbers could be defined only for the excursions of $B^H,$ and we need for its study more sophisticated methods.
\end{rem}

\vspace{7pt}

{\bf Proof of Theorem 2.} By $5)$ we have
$$
\mathbb{P}\{L(t) \neq 0,\;\; \forall\; t\in \mathbb{R}^1 \}= 1.
$$
Hence as above we can define for $L$ the winding numbers
$\nu_{\stackrel{+}{-}}^L(-\infty, t],$ 
$\nu_{\stackrel{+}{-}}^L[t, \infty),$ and besides we have
$$
\nu_{\stackrel{+}{-}}^L(-\infty, t] =  
\nu_{\stackrel{+}{-}}(0, e^t] ,\;\;\;\; 
\nu_{\stackrel{+}{-}}^L[t, \infty) = 
\nu_{\stackrel{+}{-}}[e^t, \infty). 
$$
Therefore from now on we can work with the process $L$ and will omit the index $L$ in the notation of winding numbers.
\vspace{7pt}

Let us show that 
\begin{equation}\label{nu+-}
\mathbb{P}\{|\nu_{\stackrel{+}{-}}[t, \infty)| = 
\infty, \;\;\forall t\in \mathbb{R}^1 \}= 1.
\end{equation} 
 
By  symmetry (property 3)) it is sufficient to state that 

\begin{equation}\label{nu+}
\mathbb{P}\{\nu_{+}[t, \infty) = 
\infty,\;\;\forall t\in \mathbb{R}^1 \}= 1.
\end{equation}
% So we can  simplify more the notation and replace 
%
%$\nu_{+}[t, \infty)$  by $\nu_{+}(t).$  

Using the arguments from the proof of Th. 1 we remark that the process $L$ visits infinitely often  each of four basic quadrants. It follows by continuity that at least one of 
two events $A, B$,
$$
A = \{\exists t>0,\;\; \sth \;\;\;\;\arg X(t)-\arg X(0) > \frac{\pi}{2}\},
$$
$$
B = \{\exists t>0,\;\; \sth \;\;\;\;\arg X(t)-\arg X(0) < \frac{\pi}{2}\},
$$
has probability 1. By symmetry (property 3)) 
$\mathbb{P}(A) = \mathbb{P}(B)$.
Thus, 
$$
\mathbb{P}\{\exists t>0,\;\; \sth \;\;\;\; \arg X(t)-\arg X(0)>\frac{\pi}{2}\}=1.
$$
From this follows by stationarity that 
for all $s \in \mathbb{R}^{1}$, 
$$
\mathbb{P}\{\exists t>s,\;\; \sth \;\;\;\; \arg X(t)-\arg X(s)>\frac{\pi}{2}\}=1.
$$
The set 
$$
E=\{(s,\omega)\in \mathbb{R}^{1}\times \Omega \,\;|\;\, \exists t>s, \;\; \sth \;\;\;\;\arg X(t)-\arg X(s)>\frac{\pi}{2}\}
$$
is measurable as the process 
$s\rightarrow \sup_{t>s} (\arg X(t)-\arg X(s))$ is continuous.

Based on the aforementioned and due to Fubini theorem, the set $E$ is such that $\lambda \times  \mathbb{P}(E^{\complement})=0$,
$\lambda$ being the Lebesgue measure. Therefore there exists 
$\Omega' \subset \Omega,\;\; \mathbb{P}(\Omega')=1$ such that for each 
$\omega \in \Omega'$, for almost all $s \in \mathbb{R}^{1}$, there exists $t>s$ for which $\arg X(t)-\arg X(s)>\frac{\pi}{2}$.  Take $\omega \in \Omega'$. Let us denote $E_{\omega}$ the corresponding $\omega$-section of $E$. Without loss of generality, we can suppose that for each 
$\omega \in \Omega'$, the point $0$ belongs to  $E_{\omega}$.
 As 
 $\lambda \times  \mathbb{P}(E_{\omega}^{\complement})=0$, $E_{\omega}$ is dense in $\mathbb{R}^1.$ Let 
 $u>0$ be such that 
 $\arg X(u)-\arg X(0)>\frac{\pi}{2}$. By continuity, 
 $\arg X(t)-\arg X(0)>\frac{\pi}{2}$ for all $t$ in a sufficiently small neighborhood of $u$ and therefore, there exists $t_{1} \in E_{\omega}$ for which 
 $\arg X(t_{1})-\arg X(0)>\frac{\pi}{2}$. Repeating this reasoning, we can build an increasing sequence $(t_{n})$ such that $t_{1}=0$  and $t_{n} \in E_{\omega}$. Since 
 for each  $n,\;\; \arg X(t_{n})-\arg X(t_{n-1})>\frac{\pi}{2}$, we get 
 $\sup_{t>0} (\arg X(t)-\arg X(0))=+\infty$.
 
 Thus, it is proved that for each $t$
\begin{equation}\label{wind1}
\mathbb{P}\{\nu_+[t,\infty) = +\infty\} = 1.
\end{equation} 
Now to show that
$$
\mathbb{P}\{\nu_+[t,\infty) = +\infty,\;\;\;\forall t 
\in \mathbb{R}^1\} = 1
$$ 
it is sufficient to remark that for each $\omega$ from $\Omega'$ the 
$\omega$-section $E_\omega =\mathbb{R}^1. $
Indeed, supposing that there exists $u\in E_\omega^\circ\;\;$ 
 we should have 
$$
\arg X(s)-\arg X(u)\leq \frac{\pi}{2}
$$
for each $s>t$, but that is in contradiction with the existence of
$t\in E_\omega, \; t>u,$ for which  (\ref{wind1}) holds.
Thus (\ref{nu+}) is proved.
Applying the previous reasonings to the process 
$
\{L(-t),\; t \in \mathbb{R}^1\}, 
$
we prove the remaining equalities of (\ref{winding}).
\CQFD

\vspace{10pt}

%\newpage

\end{document}